\documentclass[12pt]{amsart}

\usepackage{amssymb}

\usepackage{enumerate}

\usepackage{graphicx}

\makeatletter
\@namedef{subjclassname@2010}{%
  \textup{2010} Mathematics Subject Classification}
\makeatother

\newtheorem{theorem}{ Main Theorem}[section]
\newtheorem{thm}{Theorem}[section]

\newtheorem{lem}[thm]{Lemma}

\theoremstyle{definition}

\newtheorem{rem}[thm]{Remark}
\newtheorem{exa}{Example}
\newtheorem{pro}[thm]{Proof}

\numberwithin{equation}{section}

\begin{document}


\baselineskip=17pt


\title[On the Diophantine equations  $X^3+Y^3+Z^3+aU^k= \sum_{i=0}^na_i U_{i} ^{t_{i}}$ ; $k=3,4$ ]{On the Diophantine equations $X^3+Y^3+Z^3+aU^k= \sum_{i=0}^na_i U_{i} ^{t_{i}}$ ; $k=3,4$  }

\author[F. Izadi]{Farzali Izadi}
\address{Farzali Izadi \\
Department of Mathematics \\ Faculty of Science \\ Urmia University \\ Urmia 165-57153, Iran}
\email{f.izadi@urmia.ac.ir}

\author[M. Baghalaghdam]{Mehdi Baghalaghdam}
\address{Mehdi Baghalaghdam \\
Department of Mathematics\\ Faculty of Science \\ Azarbaijan Shahid Madani University\\Tabriz 53751-71379, Iran}
\email{mehdi.baghalaghdam@azaruniv.edu}

\date{}

\begin{abstract}
In this paper, elliptic curves theory is used for solving the Diophantine equations
$X^3+Y^3+Z^3+aU^k= \sum_{i=0}^na_i U_{i} ^{t_{i}}$,
 $k=3,4$, where $n, t_i\in \mathbb{N}\cup\{0\}$, and $a\neq 0$, $a_i$, are fixed arbitrary rational numbers. 
We try to transform each case of  the above Diophantine equations to a cubic elliptic curve of positive rank, then get infinitely many integer solutions for each case. We also solve these Diophantine equations for some values of $n$, $a$, $a_i$, $t_i$, and obtain infinitely many solutions for each case, and show among the other things that how  sums of four, five, or more cubics can be written as sums of four, five, or more biquadrates as well as sums of 5th powers, 6th powers and so on.
\end{abstract}



\maketitle

\section{Introduction}
\noindent The authors in two different papers, used  elliptic curves to solve two Diophantine equations
\begin{equation}
X^4+Y^4= 2U^4+ \sum_{i=1}^nT_i U_{i} ^{\alpha_{i}},
\end{equation}
 \noindent  where $n, \alpha_i\in \mathbb{N}$, and $T_i$, are appropriate rational numbers, and\\
 \begin{equation}
 \sum_{i=1}^n a_ix_{i} ^6+\sum_{i=1}^m b_iy_{i} ^3= \sum_{i=1}^na_iX_{i}^6\pm\sum_{i=1}^m b_iY_{i} ^3,
 \end{equation}   
 \noindent where $n, m \in \mathbb{N}$ and $a_i, b_i\neq 0\in\mathbb{Q}$. (see \cite{1}, \cite{2})\\ 

\noindent In this paper, we are interested in the study of the Diophantine equations: 
\begin{equation}\label{a} 
X^3+Y^3+Z^3+aU^k =\sum_{i=0}^na_i U_{i} ^{t_{i}};  k=3,4
\end{equation}
\noindent where  $n, t_i\in \mathbb{N}\cup\{0\}$, and $a\neq 0$, $a_i$, are fixed arbitrary rational numbers.\\
\noindent This shows that how  sums of four, five, or more cubics can be written as sums of four, five, or more biquadrates as well as sums of 5th powers, 6th powers and so on.\\

\noindent We conclude this introduction with a standard fact which is needed in section $3$. (see \cite{4})
\begin{lem} Let K be a field of characteristic not equal to $2$. Consider the equation

\noindent $v^2=au^4+bu^3+cu^2+du+q^2$, with $a$, $b$, $c$, $d$ $\in K$.
\\

\noindent Let $x=\frac{2q(v+q)+du}{u^2}$, $y=\frac{4q^2(v+q)+2q(du+cu^2)-(\frac{d^2u^2}{2q})}{u^3}$.
\\

\noindent Define
\noindent $a_1=\frac{d}{q}$, $a_2=c-(\frac{d^2}{4q^2})$, $a_3=2qb$, $a_4=-4q^2a$, $a_6=a_2a_4$.
\\

\noindent Then $y^2+a_1xy+a_3y=x^3+a_2x^2+a_4x+a_6$.
\\

\noindent The inverse transformation is
\noindent $u=\frac{2q(x+c)-(\frac{d^2}{2q})}{y}$, $v=-q+\frac{u(ux-d)}{2q}$.
\\

\noindent The point $(u, v)=(0,q)$ corresponds to the point $(x, y)=\infty$ and

\noindent $(u, v)=(0,-q)$ corresponds to $(x, y)=(-a_2,a_1a_2-a_3)$.
\end{lem}
\section{The Diophantine equation (DE) $X^3+Y^3+Z^3+aU^3 =\sum_{i=0}^na_i U_{i} ^{t_{i}}$}
\begin{theorem}Consider the DE \eqref{a} for the case $k=3$,
where $n, t_i\in \mathbb{N}\cup\{0\}$, $a\neq 0$, $a_i$, are fixed arbitrary rational numbers.\\

\noindent We try to transform this DE to a cubic elliptic curve of positive rank.\\
\noindent  Let $Y^2=X^3+FX^2+GX+H$, be an elliptic curve in which the coefficients $F$, $G$, and $H$, are all functions
of $a$, $a_i$, $t_i$ and the other rational parameters $P_i$, $t$, $s$, yet to be found later. If the elliptic curve has positive rank, depending on the values of $P_i$, $s$, the DE \eqref{a} has infinitely many  integer solutions. 
\end{theorem}

\begin{pro} Firstly, it is clear that if we find rational solutions for each case of the Diophantine equations  \eqref{a}, then by canceling the denominators of $X$, $Y$, $Z$, $U$, $U_i$, and  by multiplying  both sides of these Diophantine equations by  appropriate number $M$, we may obtain integer solutions for each case.
\\

\noindent Let $X=-Z+t$, $Y=-Z-t$,
where $Z$, $t$ are rational variables. By substituting these variables in the DE \eqref{a}, and some simplifications, we get:
 \begin{equation}\label{8}
t^2=(\frac{a}{6Z})U^3+(\frac{-Z^2}{6}-\frac{\sum_{i=0}^na_i U_{i} ^{t_{i}}}{6Z}).
\end{equation}

\noindent By multiplying  both sides of the Eq. \eqref{8}, by $(\frac{a}{6Z})^2$, and letting 
\begin{equation}\label{9}
Y'=t.(\frac{-a}{6Z}),
\end{equation}
and 
\begin{equation}\label{80}
 X'=U.(\frac{a}{6Z}),
\end{equation}
\noindent we get:
\begin{equation}\label{81}
 Y'^2=X'^3+(-\frac{a^2}{216}-(\frac{a^2}{216Z^3}).\sum_{i=0}^na_i U_{i} ^{t_{i}}).
\end{equation}

\noindent Now by choosing appropriate values for $Z$ (equal to $s$), and  $U_i$ (equal to $P_i$), such that  the rank of the elliptic curve
\eqref{81} to be positive, and by calculating $t$, $U$, $X$, $Y$, from the relations \eqref{9}, \eqref{80}, $X=-Z+t$ and $Y=-Z-t$,
 some simplifications  and canceling  the  denominators of $X, Y, Z, U, U_i$, we obtain infinitely many integer solutions for the  DE \eqref{a}. The proof is completed.
 \end{pro}

\noindent Now we work out some examples.
\\
\subsection{Application to examples} 
\begin{exa}$X^3+Y^3+Z^3+U^3=U_0^3+U_1^3+U_2^3$\\
 \noindent i.e., the sum of $4$ cubics can be written as the sum of $3$ cubics.\\
 
 \noindent The cubic elliptic curve \eqref{81} is $Y'^2=X'^3-\frac{1}{216}-(\frac{U_0^3+U_1^3+U_2^3}{216Z^3})$.\\
\noindent By taking $Z=5$, $U_0=1$, $U_1=2$, $U_2=3$, the above elliptic curve becomes\\
\noindent $Y'^2=X'^3-\frac{161}{27000}$.\\
\noindent Rank=1.\\
\noindent Generator: $P=(X',Y')=(\frac{643}{90},\frac{2578}{135})$.\\
\noindent   $(t,U)=(\frac{-5156}{9},\frac{643}{3})$.\\

\noindent Solution:\\

\noindent  $9^3+18^3+27^3+5201^3=45^3+1929^3+5111^3$.\\

\end{exa}

\begin{exa} $X^3+Y^3+Z^3+U^3=U_0^4$\\
i.e., fourth power of an integer is written as a sum of four cubics.\\
 
\noindent The cubic elliptic curve \eqref{81} is $Y'^2=X'^3-\frac{1}{216}-(\frac{U_0^4}{216Z^3})$.\\
\noindent Now letting $Z=-2$, and $U_0=7$, yields\\
\noindent $Y'^2=X'^3+\frac{2393}{1728}$.\\
\noindent Rank=2.\\
\noindent Generators: $P_1=(X',Y')=(\frac{115}{48},\frac{249}{64})$, $P_2=(X'',Y'')=(\frac{1320481}{291600},\frac{1528666129}{157464000})$.\\
 \noindent $(t,U)=(\frac{747}{16},\frac{-115}{4})$, for the point $P_1$,\\
 \noindent $(t',U')=(\frac{1528666129}{13122000},-\frac{1320481}{24300})$, for the point $P_2$.\\
 
\noindent Solutions:\\
 
\noindent $779^3+(-715)^3+(-32)^3+(-460)^3=56^4$.\\

\noindent  $ 7774550645^3+(-7512110645)^3+(-131220000)^3+(-3565298700)^3=5103000^4$.\\
\end{exa}

\begin{exa}$X^3+Y^3+Z^3+U^3=U_0^4+U_1^4+U_2^4+U_3^4$\\
 i.e., the sum of $4$ cubics is written as the sum of $4$ biquadrates.\\

\noindent The cubic elliptic curve \eqref{81} is $Y'^2=X'^3-\frac{1}{216}-(\frac{U_0^4+U_1^4+U_2^4+U_3^4}{216Z^3})$.\\
\noindent Now letting $Z=-5$, $U_0=1$, $U_1=2$, $U_2=3$, $U_3=4$, yields\\
\noindent $Y'^2=X'^3+\frac{229}{27000}$.\\
\noindent Rank=1.\\
\noindent  Generator: $P=(X',Y')=(\frac{367}{7260},\frac{14821}{159720})$.\\
\noindent $(t,U)=(\frac{14821}{5324},-\frac{367}{242})$.\\

\noindent  Solution:\\
 
\noindent $1823404^3+519156^3+(-1171280)^3+(-355256)^3=10648^4+21296^4+31944^4+42592^4$.\\
\end{exa}

\begin{exa}$X^3+Y^3+Z^3+U^3=U_0^5+U_1^5+U_2^5+U_3^5$\\
\noindent i.e., the sum of $4$ cubics is written as the sum of $4$ fifth powers.\\

\noindent The cubic elliptic curve \eqref{81} is $Y'^2=X'^3-\frac{1}{216}-(\frac{U_0^5+U_1^5+U_2^5+U_3^5}{216Z^3})$.\\
\noindent Now letting $Z=-6$, $U_0=1$, $U_1=2$, $U_2=3$, $U_3=4$, yields\\
\noindent $Y'^2=X'^3+\frac{271}{11664}$.\\
\noindent Rank=1.\\
 \noindent Generator: $P=(X',Y')=(\frac{37}{324},\frac{917}{5832})$.\\
\noindent $(t,U)=(\frac{917}{162},-\frac{37}{9})$.\\
 
 \noindent Solution:\\
 
\noindent $90672^3+2640^3+(-46656)^3+(-31968)^3=216^5+432^5+648^5+864^5$.\\
\end{exa}

\begin{exa} $X^3+Y^3+Z^3+U^3=U_0^6+U_1^6+U_2^6+U_3^6$\\
\noindent  i.e., the sum of $4$ cubics is written as the sum of $4$ sixth powers.\\

\noindent The cubic elliptic curve \eqref{81} is $Y'^2=X'^3-\frac{1}{216}-(\frac{U_0^6+U_1^6+U_2^6+U_3^6}{216Z^3})$.\\
\noindent Now letting $Z=-7$, $U_0=1$, $U_1=2$, $U_2=3$, $U_3=4$, yields\\
\noindent $Y'^2=X'^3+\frac{4547}{74088}$.\\
\noindent Rank=2.\\
 \noindent Generators:\\
\noindent  $P_1=(X',Y')=(\frac{284747}{318402},\frac{111956735}{127042398})$, \noindent $P_2=(X'',Y'')=(\frac{-162932615}{499254336},\frac{1819882008883}{11155338883584})$. \\
 \noindent $(t,U)=(\frac{111956735}{3024819},-\frac{284747}{7581})$, $(t',U')=(\frac{1819882008883}{265603306752},\frac{162932615}{11887008})$.\\
 
\noindent Solutions:\\ 

\noindent $2529478892^3+(-1724877038)^3+(-402300927)^3+(-2158667007)^3=\\
7581^6+15162^6+22743^6+30324^6$.\\

\noindent $489320985767551^3+5232371404673^3+484195324491480^3
+\\(-247276678586112)^3=$
 $5943504^6+11887008^6+17830512^6+23774016^6$.\\

\noindent If we take $Z=-7$, $U_0=U_1=3$, $U_2=U_3=4$, the above elliptic curve becomes\\
\noindent $Y'^2=X'^3-\frac{1}{216}-(\frac{3^6+3^6+4^6+4^6}{216.(-7)^3})=X^3+\frac{9307}{74088}$.\\
\noindent Rank=2.\\
\noindent Generators:\\
\noindent $P_1=(X',Y')=(\frac{1654081}{3572100},\frac{3201764021}{6751269000})$, \noindent $P_2=(X'',Y'')=(\frac{1526057}{71442},\frac{-666521155}{6751269})$.\\ 
\noindent $(t,U)=(\frac{3201764021}{160744500},-\frac{1654081}{85050})$, \noindent $(t',U')=(\frac{-1333042310}{321489},-\frac{1526057}{1701})$.\\
 
\noindent Solutions: (for the DE $X^3+Y^3+Z^3+U^3=2U_0^6+2U_2^6$):\\ 

\noindent $21634877605^3+(-10382762605)^3+(-5626057500)^3+(-15631065450)^3=\\$
\noindent $2.(85050)^6+2.(113400)^6$.\\

\noindent $(-1330791887)^3+1335292733^3+(-2250423)^3+(-288424773)^3=$\\
$2.(1701)^6+2.(2268)^6$.\\
\end{exa}
  
 \begin{exa} $X^3+Y^3+Z^3+U^3=U_0^7+U_1^7+U_2^7+U_3^7$\\
\noindent i.e., the sum of $4$ cubics is written as the sum of $4$ seventh powers.\\

\noindent The cubic elliptic curve \eqref{81} is $Y'^2=X'^3-\frac{1}{216}-(\frac{U_0^7+U_1^7+U_2^7+U_3^7}{216Z^3})$.\\
\noindent Now letting $Z=2$, $U_0=1$, $U_1=2$, $U_2=3$, $U_3=4$, yields\\
\noindent $Y'^2=X'^3-\frac{1559}{144}$.\\
\noindent Rank=1.\\
\noindent  Generator: $P=(X',Y')=(\frac{55825}{1521},\frac{-52754017}{237276})$.\\
\noindent $(t,U)=(\frac{52754017}{19773},\frac{223300}{507})$.\\

\noindent Solution:\\
 
\noindent  $365855455514133^3+(-366404379540849)^3+274462013358^3$
$+\\60441190910100^3=59319^7+118638^7+177957^7+237276^7$.\\
\end{exa}

\begin{exa} $X^3+Y^3+Z^3+U^3+U_0^3=U_1^8+U_2^8+U_3^8+U_4^8+U_5^8$\\
\noindent i.e., the sum of $5$ cubics is written as the sum of $5$ eighth powers.\\

\noindent The cubic elliptic curve \eqref{81} is $Y'^2=X'^3-\frac{1}{216}-(\frac{-U_0^3+U_1^8+U_2^8+U_3^8+U_4^8+U_5^8}{216Z^3})$.\\
\noindent Now letting $Z=2$, $U_0=1$, $U_1=\frac{2}{3}$, $U_2=2$, $U_3=\frac{1}{3}$, $U_4=3$, $U_5=\frac{4}{3}$, yields\\
\noindent $Y'^2=X'^3-\frac{14946019}{3779136}$.\\
\noindent Rank=1.\\
\noindent  Generator: $P(X',Y')=(\frac{4905677718694537513}{2497113793896353424},\frac{-7515135371775048887228789899}{3946003812454638191475172032})$.\\

\noindent $(t,U)=(\frac{7515135371775048887228789899}{328833651037886515956264336},\frac{4905677718694537513}{208092816158029452})$.\\

\noindent Solution:\\

\noindent $X=\frac{6857468069699275855316261227}{328833651037886515956264336}$,\\

\noindent $Y=\frac{-8172802673850821919141318571}{328833651037886515956264336},$\\

\noindent $U=\frac{4905677718694537513}{208092816158029452}$,\\

\noindent $Z=2$, $U_0=1$, $U_1=\frac{2}{3}$, $U_2=2$, $U_3=\frac{1}{3}$, $U_4=3$, $U_5=\frac{4}{3}$.\\
 
\end{exa}

\begin{rem}
By choosing the other points on the above elliptic curves such as 
nP ($n=3, 4$, $\cdots$,), P is one of the elliptic curves generators, we obtain infinitely many solutions for  each case of the above  Diophantine equations.
\end{rem}
\section{The DE $X^3+Y^3+Z^3+aU^4 =\sum_{i=0}^na_i U_{i} ^{t_{i}}$}
\begin{theorem} Consider the DE \eqref{a} for the case $k=4$,
where  $n, t_i\in \mathbb{N}\cup\{0\}$, and $a\neq 0$, $a_i$, are fixed arbitrary rational numbers.\\

 \noindent Let $Y^2=X^3+FX^2+GX+H$, be an elliptic curve in which the coefficients $F$, $G$, and $H$, are all functions
of $a$, $a_i$, $t_i$ and the other rational parameters $P_i$, $t$, $s$, yet to be found later. If the elliptic curve has positive rank, depending on the values of $P_i$, $s$, the DE \eqref{a} has infinitely many  integer solutions. 
\end{theorem}

\begin{pro} 
Let $X=-Z+t$, $Y=-Z-t$,
where $Z$, $t$ are rational variables. By substituting these variables in the DE \eqref{a}, and some simplifications, we get:
 \\
 
\begin{equation}\label{600}
t^2=(\frac{a}{6Z})U^4+(\frac{-Z^2}{6}-\frac{\sum_{i=0}^na_i U_{i} ^{t_{i}}}{6Z}).
\end{equation}
\\

\noindent Now by choosing appropriate values for $U_i$ (equal to $P_i$), and $Z$ (equal to $s$), such that  the rank of the  quartic elliptic curve
\eqref{600} to be positive (Then we get infinitely many rational solutions for $U$, $t$.), by calculating $U$, $t$, $X$, $Y$, $U_i$, $Z$, from the relations \eqref{600}, $X=-Z+t$, $Y=-Z-t$, $U_i=P_i$, $Z=s$, 
 after  some simplifications  and canceling  the  denominators of the values obtained for variables, we obtain infinitely many integer solutions for the  DE \eqref{a}. The proof is completed.
 \end{pro}

\begin{rem} (If in the quartic elliptic curve \eqref{600},
\begin{equation}
Q:=(\frac{-Z^2}{6}-\frac{\sum_{i=0}^na_i U_{i} ^{t_{i}}}{6Z}),
 \end{equation}
 to be square (It is done by choosing  appropriate values for $Z$, and $U_i$.), say $q^2$, we  may use the lemma $1.1$ for transforming this quartic to a cubic  elliptic curve of the form
 $y^2+a_1xy+a_3y=x^3+a_2x^2+a_4x+a_6$, where $a_i\in \mathbb{Q}$. \\
Then we solve the cubic elliptic curve just obtained of the rank$\geq1$, and get infinitely many solutions for the above DE.\\

\noindent Generally, it is not essential that $Q$ to be square, because we may transform the quartic \eqref{600} to a quartic that its  constant number is square if the rank of the quartic  \eqref{600} is positive. The only important thing is that the rank of the quartic elliptic curve \eqref{600}, to be positive for getting infinitely many nontrivial solutions for the above DE. See the example $8$.)
\end{rem}
\noindent Now let us take
\begin{equation}
 q^2=(\frac{-Z^2}{6}-\frac{\sum_{i=0}^na_i U_{i} ^{t_{i}}}{6Z}).
 \end{equation}
 
 \noindent Then for the quaratic elliptic curve \eqref{600}, we have
\begin{equation}
t^2=(\frac{a}{6Z})U^4+q^2.
\end{equation}

\noindent With the inverse transformation
\begin{equation}\label{201}
 U=\frac{2qX'}{Y'}, 
 \end{equation}
 and
\begin{equation}\label{202}
t=-q+\frac{U^2X'}{2q},
\end{equation}

\noindent the corresponding cubic  elliptic curve is
 
 \begin{equation}\label{203}
 Y'^2=X'^3+(\frac{aZ}{9}+\frac{a.(\sum_{i=0}^na_i U_{i} ^{t_{i}})}{9Z^2})X'.
 \end{equation}
\\

\noindent Then if the elliptic curve \eqref{203} has positive rank ( This is done by choosing appropriate values for $U_i$, and $Z$.), by calculating $U$, $t$, $X$, $Y$, from the relations \eqref{201}, \eqref{202}, $X=-Z+t$, and $Y=-Z-t$,
 after  some simplifications  and canceling  the  denominators of  values obtained for variables, we obtain infinitely many integer solutions for the  DE \eqref{a}. Thus we conclude that we must choose appropriate values for
$U_i$, and $Z$, such that 
 the rank of the elliptic curve \eqref{203} to be positive, then obtain infinitely many solutions for the DE \eqref{a}.
 
\subsection{Application to examples}

\begin{exa} $X^3+Y^3+Z^3=U^4+a_0U_0^{t_0}+a_1U_1^{t_1}$\\

\noindent First of all, we solve the DE 
$X^3+Y^3+Z^3=U^4+a_0U_0^{t_0}+a_1$.\\
Let $X=-Z+t$ and $Y=-Z-t$. By substituting these variables in the above DE, and some simplifications, we get

\begin{equation}\label{71}
t^2=(\frac{-1}{6Z})U^4+(\frac{-Z^2}{6}-\frac{a_0U_0^{t_0}+a_1}{6Z}).
\end{equation}
\\
\noindent Let us transform the quaratic \eqref{71} to a cubic elliptic curve.
 
\noindent Let 
\begin{equation}\label{72}
q^2=(\frac{-Z^2}{6}-\frac{a_0U_0^{t_0}+a_1}{6Z}). 
\end{equation}

\noindent Now with the inverse transformation
\begin{equation}\label{74}
 U=\frac{2qX'}{Y'}, 
 \end{equation}
 and
\begin{equation}\label{75}
t=-q+\frac{U^2X'}{2q},
\end{equation}

\noindent the quaratic \eqref{71}, maps to the cubic  elliptic curve 
 
\begin{equation}\label{77}
 Y'^2=X'^3+(\frac{2q^2}{3Z})X'.
 \end{equation}

\noindent Then if the above elliptic curve has positive rank (This is done by choosing appropriate arbitrary values for $q$, and $Z$, such as the rank of the above elliptic curve to be $\geq 1$.), we obtain infinitely many rational solutions for the DE
$X^3+Y^3+Z^3=U^4+a_0U_0^{t_0}+a_1$.\\
\noindent After multiplying both sides of the above DE by appropriate number $M$, we obtain integer solutions for the main DE
$X^3+Y^3+Z^3=U^4+a_0U_0^{t_0}+a_1U_1^{t_1}$.\\

 \noindent As an example, if  we take $Z=-2$, $q=7$, $a_0=1$, $U_0=2$, $t_0=t_1=4$, ($a_1=580$),\\ the elliptic curve \eqref{77} becomes 

\noindent $Y'^2=X'^3-\frac{49}{3}X'$.\\
\noindent Rank=1.\\
\noindent Generator: $P=(X',Y')=(\frac{-300}{361},\frac{24730}{6859})$.\\
\noindent $(t,U)=(\frac{-46590103}{6115729},\frac{-7980}{2473})$.\\

\noindent Solution:\\
 
\noindent $(-210128161627205)^3+359736726432969^3+(-74804282402882)^3=\\
48803517420^4+580.(15124197817)^4+30248395634^4$.\\
\\

\noindent By taking $Z=-1$, $q=12$, $a_0=1$, $U_0=3$, $t_0=t_1=6$, ($a_1=136$), the elliptic curve \eqref{77} becomes 
\noindent $Y'^2=X'^3-96X'$.\\
\noindent Rank=1.\\
\noindent  Generator: $P=(X',Y')=(-8,16)$.\\
 \noindent $2P=(25,-115)$.\\
 \noindent $(t,U)=(\frac{8652}{529},\frac{-120}{23})$, (for the point 2P).\\
 
\noindent  Solution:\\
 
\noindent  $4856749^3+(-4297067)^3+(-279841)^3=
63480^4+1587^6+136.(529)^6$.\\
\\

\noindent Letting $Z=-1$, $q=9$, $a_0=1$, $U_0=2$, $t_0=t_1=8$, ($a_1=231$), in the elliptic curve \eqref{77}, yields\\
\noindent $Y'^2=X'^3-54X'$.\\
\noindent Rank=1.\\
\noindent Generator: $P=(X',Y')=(-2,10)$.\\
\noindent $(t,U)=(\frac{-261}{25},\frac{-18}{5})$.\\

\noindent Solution:\\
 
\noindent $(-3687500)^3+4468750^3+(-390625)^3=
56250^4+250^8+231.(125)^8$.
\end{exa}

\begin{exa}$Y_1^3+Y_2^3+Y_3^3=X_1^5+X_2^5+X_3^5$ ($a=0$)\\
\noindent i.e., the sum of $3$ cubics can be written as the sum of $3$ fifth powers.\\

\noindent We solve this DE with another new method.\\
\noindent Let: $Y_1=t+v$, $Y_2=t-v$, $Y_3=\beta t$, $X_1=t+x_1$, $X_2=t-x_1$, $X_3=\alpha t$. \\

\noindent By substituting these variables in the above DE, we get 
\begin{equation}
(2t^3+6tv^2)+\beta^3t^3=(2t^5+10x_1^4t+20x_1^2t^3)+\alpha^5t^5.
\end{equation}
\noindent Then after some simplifications and clearing the case of $t=0$, we obtain
\begin{equation}\label{25}
v^2=\frac{2+\alpha^5}{6}t^4+\frac{20x_1^2-2-\beta^3}{6}t^2+(\frac{5}{3})x_1^4.
\end{equation}
 
\noindent  Note that $(\frac{5}{3})x_1^4$, is not a square. Then we may not use from the lemma $1.1$ for transforming this quartic to a cubic elliptic curve, but we do this work by another method. Let us take $x_1=1$, $\alpha=\beta=2$. Then the quartic \eqref{25} becomes
 
 \begin{equation}\label{26} 
 v^2=\frac{17}{3}t^4+\frac{5}{3}t^2+\frac{5}{3}.
\end{equation}

\noindent By searching, we see that the above quartic has two rational points\\ $P_1=(1,3)$, and $P_2=(7,117)$, among others. Let us put $T=t-1$. Then we get

\begin{equation}\label{27} 
 v^2=\frac{17}{3}T^4+\frac{68}{3}T^3+\frac{107}{3}T^2+26T+9.
\end{equation}

\noindent Now with the inverse transformation
\begin{equation}\label{28}
 T=\frac{6(X+\frac{107}{3})-\frac{26^2}{6}}{Y},
 \end{equation}
 and
\begin{equation}\label{29}
v=-3+\frac{T(XT-26)}{6},
\end{equation}

\noindent the quaratic \eqref{27}, maps to the cubic  elliptic curve 
 
 \begin{equation}\label{30}
 Y^2+\frac{26}{3}XY+136Y=X^3+\frac{152}{9}X^2-204X-\frac{10336}{3}.
 \end{equation}

\noindent Rank=2.\\
 \noindent Generators:
\noindent  $P_1=(X',Y')=(\frac{-152}{9},\frac{280}{27})$, $P_2=(X'',Y'')=(\frac{-44}{3},\frac{20}{9})$.\\
 
 \noindent To square the left hand of \eqref{30}, let us put $M=Y+\frac{13}{3}X+68$. Then the cubic elliptic curve \eqref{30} transforms to the Weierstrass form
 
\begin{equation}\label{31}
 M^2=X^3+\frac{107}{3}X^2+\frac{1156}{3}X+\frac{3536}{3}.
\end{equation} 

\noindent Rank=2.\\
 \noindent Generators: $G_1=(X',M')=(\frac{-44}{3}, \frac{20}{3})$, $G_2=(X'',M'')=(\frac{-152}{9},\frac{140}{27})$.\\
 
\noindent  Thus we conclude that we could transform the main quartic \eqref{26} to the cubic elliptic curve \eqref{31} of the rank equal to $2$.
  \\
  
\noindent  Because of this, the  above cubic elliptic curve has infinitely many rational points and we may obtain infinitely many solutions for the above DE too. \\
 Since $G_1=(X',M')=(\frac{-44}{3},\frac{20}{3})$, we get $(t,v)=(7,-117)$, that is on the \eqref{26}, by calculating $X_i$, $Y_i$,  from the above relations and after  some simplifications  and canceling  the  denominators of $X_i$, $Y_i$, we obtain a solution for the above DE: \\
 
\noindent   $(-110)^3+124^3+14^3=8^5+6^5+14^5$.\\

\noindent It is iteresting to see that $(-110)+124+14=8+6+14$, too.
\\

\noindent  Also we have $2G_2=(\frac{373}{36},\frac{-21721}{216})$.
 \\
 
\noindent  By using this new point $2G_2=(\frac{373}{36},\frac{-21721}{216})$, we get $(t,v)=(\frac{11}{47},\frac{2943}{2209})$.\\

\noindent Solution:\\
 
 \noindent  $359227580^3+(-251874598)^3+107352982^3=128122^5+(-79524)^5+48598^5$.\\
\end{exa}
\noindent  The Sage software has been used for calculating  the rank of  the elliptic curves. (see \cite{3})
\\


\begin{thebibliography}{HD}
\bibitem{1} F. Izadi and M. Baghalaghdam, On the Diophantine equation  $X^4+Y^4=2U^4+ \sum_{i=1}^nT_i U_{i} ^{\alpha_{i}}$, submitted, $(2016)$.
\bibitem{2} F. Izadi and M. Baghalaghdam, On the six-three power Diophantine equations, submitted, $(2016)$. 
\bibitem{3} SAGE software, available from http://sagemath.org.
\bibitem{4}  L. C. Washington, Elliptic Curves: Number Theory and Cryptography, Chapman-Hall, $(2008)$.
\end{thebibliography}
\end{document}